\documentclass{amsart}

\usepackage{amsgen,amsmath,amstext,amsbsy,amsopn,amsfonts,amssymb, amsthm}

\newtheorem{lemma}{Lemma}
\newtheorem{theorem}{Theorem}

\def\gp#1{\langle #1 \rangle}
\def\m1{^{-1}}
\author{V.~Bovdi}

\title[hyperbolic units group]{Group rings whose group of units is  hyperbolic}
\dedicatory{Dedicated to Professor L.G.~Kov\'acs on his 75th birthday}

\address{Institute of Mathematics, U. of Debrecen,
H-4010 Debrecen, P.O.B. 12, Hungary}

\email{vbovdi@math.unideb.hu}

\thanks{The research was supported by OTKA  No.K68383}



\begin{document}

\begin{abstract}
We study those group rings  whose  group of units is  hyperbolic.
\end{abstract}
\maketitle

Let $(X, \rho)$ be a metric space with  metric $\rho$. For  any $a,b,c \in  X$,
the Gromov product $\gp{b,c}_a$ of $b$ and $c$  with respect to $a\in X$  is defined as
\[
\gp{b,c}_a =\frac{1}{2}\Big(\rho (b,a)+\rho(c,a)-\rho(b,c)\Big).\textstyle
\]
The metric space is called $\delta$-hyperbolic ($\delta\geq  0$) if
\[
\gp{a,b}_d\geq \min\big\{ \gp{a,c}_d,\gp{b,c}_d\big\}-\delta \qquad(a,b,c,d\in X).
\]
Let $G$ be a finitely generated group and let $S$ be a finite set of generators of $G$. The Cayley graph
$\mathfrak{C}(G,S)$ of the group $G$ with respect to the set  $S$ is the metric
graph whose vertices are in one-to-one correspondence with the elements of $G$.
An edge  (labeled by $s$)   joins  $g$ to $gs$ for some $g\in G$ and  $s\in S$.
The group $G$ is called  {\it hyperbolic} (see \cite{Gromov}) if its Cayley graph
$\mathfrak{C}(G,S)$ is a $\delta$-hyperbolic metric space for some $\delta\geq 0$.
It is well known (see \cite{Gromov}) that this definition does not depend on the
choice of the generating set $S$.

The  following   question can be posed: When is  the group of  units $U(KG)$ of the group
ring $KG$ of a group $G$ over the commutative ring $K$ with unity  hyperbolic.
This problem was solved in \cite{Iwaki_Juriaans, Juriaans_Passi_Prasad, Juriaans_Passi_Souza_Filho}
for several particular cases.

Here we give a more general  characterization.

\begin{theorem} \label{T:1}
Let  $G$ be a group, such that the  torsion part  $t(G)\not=\{1\}$.
Let $K$ be a commutative ring  of $char(K)=0$ with unity.
If the  group of units   $U(KG)$ of the  group ring  $KG$ is hyperbolic, then  one of the following conditions holds:
\begin{itemize}
\item[(i)] $G\in \{C_5,C_8, C_{12}\}$ or $G$ is     finite  abelian of $exp(G)\in \{2,3,4,6\}$;
\item[(ii)] $G$  is a Hamiltonian $2$-group;
\item[(iii)]  $G\in\{ H_{3,2},  H_{3,4},   H_{4,2}, H_{4,4}\}$,  where $H_{s,n}=\gp{a,b\mid a^s=b^n=1, a^b=a\m1 }$;
\item[(iv)] $G=t(G)\rtimes \gp{\xi}$, where  $t(G)$ is either a finite  Hamiltonian $2$-group or
    a finite  abelian group of $exp(t(G))\in \{2,3,4,6\}$ and $\gp{\xi}\cong C_{\infty}$. Moreover, if $t(G)$ is abelian, then
    conjugation by $\xi$ either inverts all elements from $t(G)$  or leave them  fixed.
\end{itemize}
\end{theorem}

\begin{theorem} \label{T:2}
Let $KG$ be the    group algebra of a group $G$  over a field $K$ of positive characteristic,
such that the torsion part $t(G)\not=\{1\}$.  The group  of units $U(KG)$ is hyperbolic if and only if when
$K$ is a finite field  and   $G$ is a finite group.
\end{theorem}

\noindent
{\it \underline{Notation}}.
If $H\leq G$ is a subgroup of a group $G$, then  we denote by $\mathfrak{R}_l(G/H)$
the left   transversal of $G$ by  $H$.
Denote by $\Delta(G)=\{\; g\in G \mid [G:C_G(g)]<\infty\; \}$
the $FC$-center of $G$, where $C_G(g)$ is the centralizer of $g$ in $G$. The set   $t(G)$  of the torsion elements of $G$
is called the  torsion part of  $G$.
If $M\subset G$, then  we denote
the normalizer of $M$ in $G$ by
$\mathfrak{N}_{G}(M)=\{x\in G\mid xM=Mx\}$.
By $\mathfrak{Ann}_l(x)$ we denote the left annihilator of $x$ in $KG$. If $x=\sum_{g\in G}\alpha_gg\in KG$,
then  $supp(x)$ denotes  the set $\{g\in G \mid \alpha_g\not=0\}$.
If $|w|$ is   the order of $w\in t(G)$, then we put $\widehat{w}=\sum_{i=1}^{|w|}w^i\in KG$.
We denote by $C_\infty^2$ the direct product of two  cyclic groups of infinite orders.
\smallskip

We shall use freely the following properties of  hyperbolic groups:
\begin{lemma} \label{L:1}
(\cite{Bridson_Haefliger}) If  $G$ is  a hyperbolic group, then:
\begin{itemize}
\item[(i)] $C_\infty^2$ does not embed as a subgroup of $G$;
\item[(ii)] if $g\in G$  has infinite order, then $[C_G(g) : \gp{g}]$ is finite;
\item[(iii)] torsion subgroups of $G$ are finite of bounded order.
\end{itemize}
\end{lemma}
First of all we shall need the following lemmas.
\begin{lemma}\label{L:2}
If   $w\in t(G)$, then $\mathfrak{Ann}_l(\widehat{w})=KG(w-1)$, where $\widehat{w}=\sum_{i=1}^{|w|}w^i\in KG$.
\end{lemma}
\begin{proof}
Let $z\in \mathfrak{Ann}_l(\widehat{w})$ and $z=\sum_{i}z_iw_i$, where $w_i=\sum_{j}\alpha_{ij}w^j\in K\gp{w}$ and
$z_i\in \mathfrak{R}_l(G/\gp{w})$.
If $\chi$  is the augmentation map,  then
\[
0=z\widehat{w}=\sum_{i}z_iw_i\widehat{w}=\sum_{i}z_i\chi(w_i)\widehat{w},
\]
so $\chi(w_i)=\sum_{j}\alpha_{ij}=0$ for all $i$,  and
\[
\begin{split}
z=\sum_{i}z_i\sum_{j}\alpha_{ij}w^j=\sum_{i}z_i\Big(\sum_{j}&\alpha_{ij}(w^j-1)+\\
&+\sum_{j}\alpha_{ij}\Big)=\sum_{i,j}z_i\alpha_{ij}(w^j-1).
\end{split}
\]
Using the equality $xy-1=(x-1)(y-1)+(x-1)+(y-1)$, where  $x,y\in G$,  we obtain that   $\mathfrak{Ann}_l(\widehat{w})\subseteq KG(w-1)$.
The opposite inclusion  is trivial.
\end{proof}

\begin{lemma}\label{L:3}
Let $G$ be a group which has at least one non-normal finite subgroup $\gp{w}$.
Put $t_w=(w-1)(\alpha g+\beta h)\in KG$, where  $\alpha,\beta\in K\setminus\{0\}$ and  $g,h\in G\setminus \mathfrak{N}_{G}(\gp{w})$.
The product  $t_w\widehat{w}=0$ if and only if one of the following conditions holds:
\begin{itemize}
\item[(i)] $g\in h\gp{w}$\quad  and\quad   $\alpha+\beta=0$;
\item[(ii)]  $wg\in h\gp{w}$, \quad $wh\in g\gp{w}$\quad   and \quad $\alpha=\beta$.
\end{itemize}
\end{lemma}

\begin{proof}
Let $t_w\widehat{w}=0$.
Since $t_w\in \mathfrak{Ann}_l(\widehat{w})$ and $|supp(t_w)|=4$,  by
Lemma \ref{L:2} we have
\begin{equation}\label{E:1}
t_w=(g_1w_1+g_2w_2)(w-1),
\end{equation}
where  $g_1,g_2\in \mathfrak{R}_l(G/\gp{w})$, $w_1,w_2\in K\gp{w}$ and
$|supp(w_i(w-1))|\in\{0, 2, 4\}$. Assume that $g\in g_1\gp{w}$. If $|supp(w_1(w-1))|=4$, then
$wg\in g_1\gp{w}$ and $g\in \mathfrak{N}_{G}(\gp{w})$,
which is impossible.  Hence $|supp(g_i(w_i-1)(w-1))|=2$. From (\ref{E:1}) we have
\[
\alpha(wg)+\beta(wh) -\alpha g -\beta h= g_1w_1w +g_2w_2w-g_1w_1-g_2w_2.
\]
Clearly   $g, h, wg, wh \in \{ g_1\gp{w}, g_2\gp{w}\}$, so we need to consider only the following possible cases:

Case 1. Let $g, h \in g_1\gp{w}$ and $ wg,wh \in g_2\gp{w}$. Then $g=hw^k$, for some $k$ and
$t_w\widehat{w}=(\alpha+\beta)(w-1)h\widehat{w}=0$. Therefore $\alpha+\beta=0$ and (i) holds.

Case 2. Let $g, wh \in g_1\gp{w}$ and $ wg,h \in g_2\gp{w}$. Then $wh=gw^l$ and $wg=hw^m$, for some $l,m$. This yields that
\[
t_w\widehat{w}=(\alpha-\beta)(h-g)\widehat{w}=0.
\]
If $\alpha-\beta\not=0$, then $h\widehat{w}=g\widehat{w}$, so from  $wh\in g\widehat{w}$ we get that $g\in \mathfrak{N}_{G}(\gp{w})$,  a contradiction.
Consequently $\alpha=\beta$  and (ii) holds.

The converse statement is obvious.
\end{proof}

\begin{lemma}\label{L:4}
Let $G$ be a group which has at least one non-normal finite subgroup $\gp{w}$.  Let   $\mathfrak{x}_w(g)=1+(w-1)g\widehat{w}$\; and\; $\mathfrak{x}_w(h)=1+ (w-1)h\widehat{w}$, where  $g,h\in G\setminus \mathfrak{N}_{G}(\gp{w})$. If $char(K)=0$ and
$\mathfrak{x}_w(g)\not\in \{\mathfrak{x}_w(h), \mathfrak{x}_w(h)\m1\}$,
then  one of the following condition holds:
\begin{itemize}
\item[(i)]  $wg\not\in h\gp{w}$ and  $g\not\in h\gp{w}$;\qquad
(ii) $wh\not\in g\gp{w}$ and  $g\not\in h\gp{w}$.
\end{itemize}
Moreover, if either condition (i) or (ii) holds, then    $\gp{\mathfrak{x}_w(g), \mathfrak{x}_w(h)}\supseteq  C_\infty^2$.
\end{lemma}

\begin{proof} Let $\mathfrak{x}_w(g)\not\in \{\mathfrak{x}_w(h), \mathfrak{x}_w(h)\m1\}$. Then  $(w-1)(g\pm h)\widehat{w}\not=0$ and  by Lemma \ref{L:3} we obtain    either  (i) or (ii).
If $\mathfrak{x}_w(g)^\alpha=\mathfrak{x}_w(h)^\beta$ for some $\alpha,\beta\in \mathbb{Z}\setminus \{0\}$,  then
\[
\big(\alpha(wg)-\beta (wh)-\alpha g+\beta h\big)\widehat{w}=0,
\]
which leads  us  to a contradiction. Consequently, $\mathfrak{x}_w(g), \mathfrak{x}_w(h)$ are different commuting units of infinite order, such that
$\mathfrak{x}_w(g)^\alpha\not=\mathfrak{x}_w(h)^\beta$  for all  different  $\alpha, \beta \in \mathbb{Z}\setminus \{0\}$,\quad
 so $\gp{\mathfrak{x}_w(g), \mathfrak{x}_w(h)}\supseteq  C_\infty^2$.
\end{proof}

\begin{lemma}\label{L:5} (see \cite{Artamonov_Bovdi}, 2.2)
If $G$ is a finite abelian  group, then $U(\mathbb{Z}G)= \pm G\times S_G$ and  the rank $\mathfrak{r}(S_G)$
of the  free abelian subgroup $S_G$ is   $\frac{1}{2}(|G|+1+t_2-2l)$, where
$t_2$ and $l$ are the cardinalities   of the involutions  and the  cyclic subgroups of $G$, respectively.
\end{lemma}

\begin{lemma}\label{L:6}
Let $G$ be a group which has at least one non-normal finite subgroup,
and assume that the order of each finite cyclic subgroup lies in
$\{2,3,4,5,6,8,12\}$.
Assume also that for each non-normal finite cyclic subgroup $\gp{w}$
and for all $g$, $h$ outside the normalizer $N$ of $\gp{w}$,
either $g\gp{w}=h\gp{w}$ or both of the following conditions hold:
   \begin{equation}\label{E:2}
    wg\in h\gp{w}\qquad \quad   and   \qquad   \quad     wh\in g\gp{w}.
   \end{equation}
Then $G\in\{ H_{3,2},  H_{3,4},   H_{3, 8},  H_{4,2}, H_{4,4}, Q_{16}\}$, where
\[
Q_{16}=\gp{a,b\mid a^8=1, b^2=a^4, a^b=a\m1}.
\]

\end{lemma}

\begin{proof}
Let $G$, $\gp{w}$, $N$, $h$ be as in the hypotheses:
then each $g$ outside $N$ lies either in $h\gp{w}$ or in $w\m1h\gp{w}$.
This shows that the set-complement $G\setminus N$,
which is the union of the non-trivial cosets of $N$,
is covered by two cosets of $\gp{w}$.
It follows that $N$ has at most two non-trivial cosets: $|G:N|\leq3$.
When this index is 3,
the two non-trivial cosets of $N$ are covered by two cosets of $\gp{w}$,
so we must have $\gp{w}=N$.
Otherwise $N$ is normal, so $\gp{w}<N$;
the unique non-trivial coset of $N$ being covered by two cosets of $\gp{w}$ then
gives that in fact $|N:\gp{w}|=2$, $|G:\gp{w}|=4$.
Since $\gp{w}$ is finite, so is $G$; indeed, by our hypothesis on the possible
orders of $w$, we have $|G|\leq 48$. All groups so small have been well understood
for a long time, and it should not be hard to find a clever way of exploiting that
understanding here. The reader will be able to choose between several ways of completing
this proof. The method preferred by the author is to take advantage of the computer
algebra system   GAP \cite{GAP}, which contains databases of small groups and can
quickly search through them, confirming the claim of the lemma. (Of course it also
shows that all six groups in the conclusion satisfy the hypotheses.)
\end{proof}

\begin{proof}[\underline{Proof of the Theorem \ref{T:1}}]
Let $G$ be a group which  has at least one non-normal finite subgroup $\gp{w}$ and let $g,h\in G\setminus \mathfrak{N}_{G}(\gp{w})$.
By Lemma \ref{L:5},  $|w|\in  \{2,3,4,5,6,8,12\}$.

\smallskip

We divide the problem into  the following two cases:

Case 1. either $wg\not \in h\gp{w}$ or $wh\not \in g\gp{w}$; \quad Case 2. $wg\in h\gp{w}$  and     $wh\in g\gp{w}$.

We consider each case separately.

Case 1. Let $wg\not \in h\gp{w}$ or $wh\not \in g\gp{w}$. If  $g\not\in h\gp{w}$, then put
$\mathfrak{x}_w(g)=1+(w-1)g\widehat{w}$ and   $\mathfrak{x}_w(h)=1+ (w-1)h\widehat{w}$. By Lemma \ref{L:4},
$U(KG)\supset \gp{\mathfrak{x}_w(g),\mathfrak{x}_w(h)} \supseteq C_\infty^2$,
a contradiction to Lemma \ref{L:1}(i). Therefore $\gp{w}$ is  normal in $G$, a contradiction.

If  for a given $h$, each element $g$  in the set complementer $G\setminus N$, where $N=\mathfrak{N}_{G}(\gp{w})$,  satisfy   $g\gp{w}= h\gp{w}$, then $G=N\cup aN$ for some $a\in G$ and   $N=\gp{w}$. Consequently $|G:\gp{w}|=2$ and $\gp{w}$ is  normal in $G$, again  a contradiction.

Case 2. Let $wg\in h\gp{w}$  and     $wh\in g\gp{w}$. According to the previous case, we can assume that this condition holds   for all
non-normal finite subgroup $\gp{w}$  and for all $g,h\in G\setminus   \mathfrak{N}_{G}(\gp{w})$. Then $g\not\in h\gp{w}$. Indeed, if $g=hw^i$, then
from $wg\in h\gp{w}$ we have  $g\in \mathfrak{N}_{G}(\gp{w})$, a contradiction.

Now, by   Lemma \ref{L:6},    $G$ is either a Hamiltonian group or it is  finite and
\[
G\in\{ H_{3,2},  H_{3,4},   H_{3, 8},  H_{4,2}, H_{4,4}, Q_{16}\}.
\]
Let  $G=\gp{a,b\mid a^8=1, b^2=a^4, a^b=a\m1}\cong Q_{16}$. The elements  $w_1=1+(a\m1-a)\widehat{b}$ and $w_2=6(a+a^3)(1-a^4)+9a^2-8a^6$
are commuting units in $V(\mathbb{Z}G)$ of infinite order, because  $w_1^{-1}=1-(a\m1-a)\widehat{b}$ and $w_2^{-1}=6(a+a^3)(1-a^4)+9a^6-8a^2$. Since
$w_1^n=1+n(b-1)a\widehat{b}$, so $ab\in supp(w_1^n)$, but $ab\not\in supp(w_2^m)$ for all nonzero $n,m\in \mathbb{Z}$. This yields that
$w_1^n\not=w_2^m$, so $\gp{w_1,w_2} \supseteq C_\infty^2$, a contradiction.

Now, if  $G=\gp{a,b\mid a^3=b^8=1, a^b=a\m1 }\cong H_{3,8}$, then using the same argument for the units $w_1=1+(a\m1-a)\widehat{b}$ and $w_2=6(b+b^3)(1-b^4)+9b^2-8b^6$ we obtain that $\gp{w_1,w_2} \supseteq C_\infty^2$, a contradiction. Therefore $G\in\{ H_{3,2},  H_{3,4},  H_{4,2}, H_{4,4}\}$.

From the previous two cases follows that  if $G$ does not satisfy  condition (ii) of our theorem, then   each cyclic subgroup of  $t(G)$ is  normal in $G$, so  $t(G)$ is either  abelian or a Hamiltonian group. Moreover,    by Lemma \ref{L:1}(iii),  $t(G)$ is finite. If $t(G)$ is abelian,
then     by Lemma \ref{L:5}  we must have  $\mathfrak{r}(S_{t(G)})\leq 1$, so either $t(G)\in \{C_5,C_8, C_{12}\}$
or  $t(G)$  is  finite  abelian of $exp(t(G))\in \{2,3,4,6\}$ (see \cite{Artamonov_Bovdi}, {\bf 3.1}), because $\mathbb{Z}G\subseteq KG$.

Let  $t(G)$ be  a Hamiltonian group such that  $\gp{a,b}\times \gp{c}\subseteq t(G)$, where $\gp{a,b}\cong Q_8$ and  $|c|=p$ is odd.
If $H=\gp{a^2c}$ then,   by   Lemma \ref{L:5}, $\mathfrak{r}(S_H)=\frac{2p-4}{2}=1$  so the only possible case is $p=3$.
Using the argument of ${\bf 4}^0$ of \cite{Juriaans_Passi_Prasad},  we get a contradiction.
Therefore if $t(G)$ is nonabelian then  it is a Hamiltonian $2$-group.

By the definition of  hyperbolic groups, the group  $U(KG)$ is
finitely  generated.
Clearly   $t(G)=t(\Delta(G))$ and  by a theorem of B.H.~Neumann (see \cite{Artamonov_Bovdi}, \S4 or  \cite{Erdos}),
the $FC$-group  $\Delta(G)/t(G)$ is   torsion free abelian. Assume that
\[
\Delta(G)/t(G)\supseteq  \gp{z_1}\times \gp{z_2}\cong C_\infty^2
\]
and let $h_1,h_2$ be the pre-images of $z_1$ and  $z_2$, respectively.
By \cite{Erdos}, we can assign to the  finitely   generated $FC$-group $\gp{h_1,h_2}$
a natural number $m$, such that  $h_1^m$ and $h_2^m$ are central in $\gp{h_1,h_2}$  and
$(h_1h_2)^m =h_1^m h_2^m$, so $\gp{h_1,h_2}\supseteq  C_\infty^2$, a contradiction. Consequently
$\Delta(G)/t(G)$ is cyclic and\quad
$\Delta(G)=t(G)\rtimes \gp{\xi}$,\quad where   $\gp{\xi}\cong C_\infty$.

Let $g\in G\setminus \Delta(G)$, such that $\gp{g}\cap \gp{\xi}=\gp{1}$. Of course $g$ is of infinite order, $\xi^g\in \Delta(G)$ and
$\Delta(G)/t(G)$ is cyclic, so we  have
$\xi^g=w\xi^{\pm 1}$, where $w\in t(G)$ and $w^g\in \gp{w}$. Since $t(G)$ is finite, there  exists an $s>1$, such that  $(\xi,g^s)=1$ and
$\gp{\xi,g^s}\cong C_\infty^2$, a contradiction. Therefore  $\xi^l=g^k\in \Delta(G)$ for some $l,k> 1$,  and
$[G:C_G(\gp{g^k})]<\infty$. Since $\gp{g^k}$ is infinite, by Lemma \ref{L:1}(ii),
$[C_G(\gp{g^k}):\gp{g^k}]<\infty$, so $[G:\gp{g}]<\infty$  and $g\in \Delta(G)$, a contradiction.  Consequently
$G=\Delta(G)$.

Now let $G=t(G)\rtimes \gp{\xi}$, where $\gp{\xi}\cong C_\infty$ and each cyclic subgroup of the finite group $t(G)$ is normal
in $G$. Clearly $\gp{\xi^k}$ is a central subgroup of $G$ for some $k>1$. If
$t(G)\supseteq  \{C_5, C_8, C_{12}\}$, then (see 2.2, 2.3,  \cite{Artamonov_Bovdi}) the subring $\mathbb{Z}t(G)\subseteq KG$ contains a  unit
$u$ of infinite order and  $\gp{\xi^k, u}\supseteq  C_\infty^2$,
a contradiction.

Consequently the part (iv) of our theorem holds. \end{proof}

\begin{lemma}\label{L:7} (\cite{Bovdi}, Proposition 2.7, p.9)
Let $H$ be  a subgroup of a group $G$. The left  annihilator  $L$
in $KG$ of the right  ideal $\mathfrak{I}_r(H)=\gp{h-1\mid h\in H,
h\not=1}$ is different  from zero if and only if $H$ is finite.
If $H$ is finite, then $L=KG(\Sigma_{h\in H}h )$.
\end{lemma}

\begin{proof}[\underline{Proof of the Theorem \ref{T:2}}]
Let $K$ be a field of $char(K)=p$.
First  assume  that $K$ contains  a transcendental element $\xi$.
If $g\in G$ is of prime order  $q$ with $(q,p)=1$, then $e=\frac{\widehat{g}}{|g|}$ and  $f=1-e$ are idempotents of $KG$,
and\quad  $\gp{\xi e+f, e+\xi f}\cong C_\infty^2$,\quad  a contradiction.
Furthermore, if  $K$ is an infinite  field such that each element of $K$ is algebraic,
then for the prime $q$  and for the set  $\{ n\mid n\in \mathbb{N}\}$ there exist infinity many different
$\beta(n)\in \mathbb{N}$  with the following property:
\[
p^n-1=q^{\alpha(n)}\beta(n)\qquad  \text{and}\qquad (q,\beta(n))=1.
\]
Since  $\mathfrak{T}=\{\alpha \in K\mid \alpha^{\beta(n)}=1,\;  n\geq 1\}$ is infinite, we obtain
$\gp{e+\alpha f\mid  \alpha\in \mathfrak{T} }$ is an infinite  torsion group, a contradiction
to Lemma \ref{L:1}(iii). Thus $K$ is finite.

Let  $w\in G$ such that $|w|=p$ and let $K$ be still  infinite. Clearly  $(\widehat{w})^p=p\widehat{w}=0$ and
 $1+\alpha \widehat{w}$, were  $\alpha\in K$,  is a unit of order $p$. The abelian  group
$\gp{1+\alpha \widehat{w}\mid \alpha \in K}$ is infinite, a contradiction  to Lemma \ref{L:1}(iii).

Consequently, we have proven that $K$ is always a finite field.

Let $G$ be a group which  has at least one non-normal finite subgroup $\gp{w}$ and  assume  that  $G\not=t(G)$. Put
\[
A=\mathfrak{R}_l(G/\mathfrak{N}_G(\gp{w}))=\{g_i\mid i\geq 1\}.
\]
First suppose that   $A$ is an infinite set.
Clearly, if $i\not=j$, then   $g_i\not\in g_j\gp{w}$. Furthermore    for each  $i\geq 1$  exists at most  one $j$ such that
$wg_i\in g_j\gp{w}$. Using $A$,  we can construct inductively  a family of the sets $\{A_i\mid i\geq 1\}$ in the following way:

Put $A_1=\{a_1\}$, where $a_1=g_1$. Consider the infinite set $A\setminus A_1$. It contains at most one element $g_j$ such that $wa_1\in g_j\gp{w}$.  We delete $g_j$ (if such element exists) from the infinite set $A\setminus A_1$. Thus we obtain an infinite set, whose elements  can be relabeled
by $\{a_2,a_3,\ldots\}$. Now put $A_2=A_1\cup \{a_2\}$. Clearly $wa_1\not \in a_2\gp{w}$.

If   we have  constructed  the  sets $A_1, \ldots, A_s$, then put
$A_{s+1}=A_{s}\cup \{a_{s+1}\}$, where  $wa_{i}\not \in a_j\gp{w}$ with $1\leq i,j\leq s+1$. By the  continuation   of  this process we obtain an infinite set
$\overline{A}=\{a_i\mid i\geq 1\}=\cup_{n\geq 1}A_n\subseteq A$,
such that  $wa_i\not\in a_j\gp{w}$ and $a_i\not\in a_j\gp{w}$ for all  $i\not=j$.
By Lemma \ref{L:3},
\[
H_w=\gp{\; \mathfrak{x}_w(g_i)=1+(w-1)g_i\widehat{w}\; \mid \; g_i\in \overline{A},\quad  i\geq 1}
\]
is an infinite abelian subgroup of  exponent $p$ and,  by Lemma \ref{L:1}(iii),  it is finite,  a contradiction.  Consequently
$A$ is a finite set and $\mathfrak{N}_G(\gp{w})$ is an infinite set.

Let $G\not=\mathfrak{N}_G(\gp{w})$ and let $c\in G\setminus\mathfrak{N}_G(\gp{w})$.
Consider the  abelian group
\[
\begin{split}
M_c=\gp{\; \mathfrak{x}_w(ch_i)=1+(w-1)ch_i\widehat{w}\quad  \mid \quad  h_i\in  \mathfrak{R}_l(\mathfrak{N}_G(\gp{w})/\gp{w}),  \quad  i\geq 1\;}.
\end{split}
\]
If $\mathfrak{x}_w(ch_i)=\mathfrak{x}_w(ch_j)$, then
$(w-1)c\widehat{w}(h_ih_j^{-1}-1)=0$.
Since $(w-1)c\widehat{w}\not=0$ and $\mathfrak{R}_l(\mathfrak{N}_G(\gp{w})/\gp{w})$ is an infinite set, we get a
contradiction to Lemma \ref{L:7}. Therefore    $M_c$ is an infinite  group
of exponent $p$, a contradiction to Lemma \ref{L:1}(iii).

Therefore $G=\mathfrak{N}_G(\gp{w})$ for all $w\in t(G)$. This  yields  that   each cyclic subgroup of $t(G)$ is
normal in  $G$, so by Lemma \ref{L:1}(iii)   the group $t(G)$ is either  finite abelian or a finite Hamiltonian group.
Clearly   $t(G)=t(\Delta(G))$ and  by a theorem of B.H.~Neumann (see \cite{Artamonov_Bovdi}, \S4 or  \cite{Erdos}),
the $FC$-group  $\Delta(G)/t(G)$ is   torsion free abelian. Assume
\[
\Delta(G)/t(G)\supseteq  \gp{z_1}\times \gp{z_2}\cong C_\infty^2
\]
and let $h_1,h_2$ be the pre-images of $z_1$ and  $z_2$, respectively.
By \cite{Erdos}, we can assign to the  finitely   generated $FC$-group $\gp{h_1,h_2}$
a natural number $m$, such that  $h_1^m$ and $h_2^m$ are central in $\gp{h_1,h_2}$  and
$(h_1h_2)^m =h_1^m h_2^m$, so $\gp{h_1,h_2}\supseteq  C_\infty^2$, a contradiction. Consequently
$\Delta(G)/t(G)$ is cyclic and
\[
G=t(G)\rtimes \gp{t},\quad \text{where}\quad   \gp{t}\cong C_\infty.
\]
Let $w\in Syl_p(t(G))$. It is easy to check that $1+t^i\widehat{w}$
is a torsion unit of order $p$ and $\gp{1+t^i\widehat{w}\mid i\in \mathbb{Z}}$
is an infinite torsion group, a contradiction  to Lemma \ref{L:1}(iii).  Thus $Syl_p(t(G))=\gp{1}$ and $t(G)$ is a $p'$-group.
If  $g\in t(G)$ is  of prime order such that $(p,|g|)=1$, then $e=\frac{\widehat{g}}{|g|}$ and  $f=1-e$ are idempotents of $KG$ and  $\gp{et+f, e+ft}\cong C_\infty^2$, a contradiction.
Consequently $t(G)=\gp{1}$,  which is impossible. Therefore   $G=t(G)$ and, by Lemma \ref{L:1}(iii),  the group $G$ is  finite.

The converse case is trivial, since $U(KG)$ is finite, so it is hyperbolic.
\end{proof}

\bibliographystyle{abbrv}
\bibliography{Hyperbolic_III}

\begin{thebibliography}{1}

\bibitem{Artamonov_Bovdi}
V.~A. Artamonov and A.~A. Bovdi.
\newblock Integral group rings: groups of invertible elements and classical
  {$K$}-theory.
\newblock In {\em Algebra. {T}opology. {G}eometry, {V}ol.\ 27 ({R}ussian)},
  Itogi Nauki i Tekhniki, pages 3--43, 232. Akad. Nauk SSSR Vsesoyuz. Inst.
  Nauchn. i Tekhn. Inform., Moscow, 1989.
\newblock Translated in J. Soviet Math. {{\bf{5}}7} (1991), no. 2, 2931--2958.

\bibitem{Bovdi}
A.~A. Bovdi.
\newblock Group rings. ({R}ussian).
\newblock {\em Kiev.UMK VO}, page 155, 1988.

\bibitem{Bridson_Haefliger}
M.~R. Bridson and A.~Haefliger.
\newblock {\em Metric spaces of non-positive curvature}, volume 319 of {\em
  Grundlehren der Mathematischen Wissenschaften [Fundamental Principles of
  Mathematical Sciences]}.
\newblock Springer-Verlag, Berlin, 1999.

\bibitem{Erdos}
J.~Erd{\"o}s.
\newblock The theory of groups with finite classes of conjugate elements.
\newblock {\em Acta Math. Acad. Sci. Hungar.}, 5:45--58, 1954.

\bibitem{GAP}
The GAP~Group.
\newblock {\em {GAP -- Groups, Algorithms, and Programming, Version 4.4.10}},
  2007.
\newblock \verb+(http://www.gap-system.org)+.

\bibitem{Gromov}
M.~Gromov.
\newblock Hyperbolic groups.
\newblock In {\em Essays in group theory}, volume~8 of {\em Math. Sci. Res.
  Inst. Publ.}, pages 75--263. Springer, New York, 1987.

\bibitem{Iwaki_Juriaans}
E.~Iwaki and S.~O. Juriaans.
\newblock Hypercentral unit groups and the hyperbolicity of a modular group
  algebra.
\newblock {\em Comm. Algebra}, 36(4):1336--1345, 2008.

\bibitem{Juriaans_Passi_Prasad}
S.~O. Juriaans, I.~B.~S. Passi, and D.~Prasad.
\newblock Hyperbolic unit groups.
\newblock {\em Proc. Amer. Math. Soc.}, 133(2):415--423 (electronic), 2005.

\bibitem{Juriaans_Passi_Souza_Filho}
S.~O. Juriaans, I.~B.~S. Passi, and A.~C. Souza~Filho.
\newblock Hyperbolic unit groups and quaternion algebras.
\newblock {\em Proc. Indian Acad. Sci. Math. Sci.}, 119(1):9--22, 2009.

\end{thebibliography}

\end{document}